\magnification=1130
\vsize=190mm
\hsize=135mm
\hoffset=3mm
\voffset=6mm
\null

\def\ConclLossCommon{2.1}
\def\RDeqn{2.2}

\def\ConvBoundary{3.1}
\def\PoincareIneq{3.2}
\def\approxpbm{3.3}
\def\contdep{3.4}
\def\classreg{3.5}
\def\estimut{3.6}
\def\Bernstein{3.7}

\def\propdeltaA{4.1}
\def\propdeltaB{4.2}
\def\integralvarphiA{4.3}
\def\integralvarphi{4.4}
\def\integralvarphiB{4.5}
\def\uWonep{4.6}
\def\integralvarphiD{4.7}

\def\AuxilPbmV{5.1}
\def\TangentBall{5.2}

\def\ratiowA{7.1}
\def\ratiowB{7.2}
\def\HypNoLoss{7.3}

\def\PsiElliptA{8.1}
\def\PsiElliptB{8.2}
\def\PsiEllipt{8.3}
\def\PsiLarge{8.4}
\def\TangentBallB{8.5}
\def\condVzero{8.6}

\def\Alaa{1}
\def\ABG{2}
\def\ARS{3}
\def\BaD{4}   
\def\BD{5}   
\def\CG{6}   
\def\CIL{7}   
\def\FL{8} 
\def\FTW{9} 
\def\Fri{10}
\def\GH{11}
\def\HM{12}
\def\KPZ{13}
\def\KS{14}
\def\LS{15}
\def\Lio{16}
\def\Ma{17}
\def\NST{18}
\def\PS{19}
\def\PZ{20}
\def\Qu{21}
\def\QS{22}
\def\Sou{23}
\def\SV{24}
\def\SZ{25}
\def\TT{26}

\def\eps{\varepsilon}

\def \trait (#1) (#2) (#3){\vrule width #1pt height #2pt depth #3pt}
\def \fin{\null\hfill
        \trait (0.1) (5) (0)
        \trait (5) (0.1) (0)
        \kern-5pt
        \trait (5) (5) (-4.9)
        \trait (0.1) (5) (0)
\medskip}
\null

\newtoks\hautpagegauche
\newtoks\hautpagedroite
\newtoks\paragraphecourant
\newtoks\chapitrecourant
\hautpagegauche={}
\hautpagedroite={}
\headline={\ifodd\pageno\the\hautpagedroite\else\the\hautpagegauche\fi}

\font\TenEns=msbm10
\font\SevenEns=msbm7
\font\FiveEns=msbm5
\newfam\Ensfam
\def\Ens{\fam\Ensfam\TenEns}
\textfont\Ensfam=\TenEns
\scriptfont\Ensfam=\SevenEns
\scriptscriptfont\Ensfam=\FiveEns
\def\R{{\Ens R}}

\def\eps{\varepsilon}


\font\itsmall=cmsl9

\font\eightrm=cmr9
\font\sixrm=cmr6
\font\fiverm=cmr5

\font\eighti=cmmi9
\font\sixi=cmmi6
\font\fivei=cmmi5

\font\eightsy=cmsy9
\font\sixsy=cmsy6
\font\fivesy=cmsy5

\font\eightit=cmti9
\font\eightsl=cmsl9
\font\eighttt=cmtt9

\def\eightpoint{\def\rm{\fam0\eightrm}
\textfont0=\eightrm
\scriptfont0=\sixrm
\scriptscriptfont0=\fiverm

\textfont1=\eighti
\scriptfont1=\sixi
\scriptscriptfont1=\fivei

\textfont2=\eightsy
\scriptfont2=\sixsy
\scriptscriptfont2=\fivesy

\textfont3=\tenex
\scriptfont3=\tenex
\scriptscriptfont3=\tenex

\textfont\itfam=\eightit \def\it{\fam\itfam\eightit}
\textfont\slfam=\eightsl \def\it{\fam\slfam\eightsl}
\textfont\ttfam=\eighttt \def\it{\fam\ttfam\eighttt}
}


\font\pc=cmcsc9
\font\itsmall=cmsl9

\paragraphecourant={\rmt DIFFUSIVE HAMILTON-JACOBI EQUATION}
\chapitrecourant={\rmt PORRETTA AND SOUPLET}
\footline={\ifnum\folio=1 \hfill\folio\hfill\fi}
\hautpagegauche={\tenrm\folio\hfill\the\chapitrecourant\hfill}
\hautpagedroite={\ifnum\folio=1 \hfill\else\hfill\the\paragraphecourant\hfill\tenrm\folio\fi}

\font\rmb=cmbx8 scaled 1125 \rm

\centerline{\rmb ANALYSIS OF THE LOSS OF BOUNDARY CONDITIONS}
\vskip 0.5mm
\centerline{\rmb FOR THE DIFFUSIVE HAMILTON-JACOBI EQUATION}

\vskip 5mm
\centerline{\pc Alessio PORRETTA}
\vskip 3mm
\centerline{\itsmall Universit\`a di Roma Tor Vergata}
\centerline{\itsmall Dipartimento di Matematica, Via della Ricerca Scientifica 1}
\centerline{\itsmall 00133 Roma, Italia. Email: porretta@mat.uniroma2.it}

\centerline{\itsmall }

\vskip 2mm
\centerline{\pc Philippe SOUPLET}
\vskip 3mm
\centerline{\itsmall Universit\'e Paris 13, Sorbonne Paris Cit\'e, CNRS UMR 7539}
\centerline{\itsmall Laboratoire Analyse G\'eom\'etrie et Applications}
\centerline{\itsmall 93430 Villetaneuse, France. Email: souplet@math.univ-paris13.fr}

\vskip 5mm

\baselineskip=12pt
\font\rmt=cmr9

\setbox1=\vbox{
\hsize=120mm
{\baselineskip=11pt \parindent=3mm \eightpoint \rmt
{\pc Abstract:}\  We consider the diffusive Hamilton-Jacobi equation,
with homogeneous Dirichlet conditions and regular initial data.
It is known from [Barles-DaLio, 2004] that the problem admits a unique, continuous, global viscosity solution, 
which extends the classical solution in case gradient blowup occurs. We study the question of the possible loss of boundary conditions
after gradient blowup, which seems to have remained an open problem till now.

Somewhat surprisingly, our results show that the issue strongly depends on the initial data
and reveal a rather rich variety of phenomena.
For any smooth bounded domain, we construct initial data such that the loss of boundary conditions occurs everywhere on the boundary,
as well as initial data for which no loss of boundary conditions occurs in spite of gradient blowup.
Actually, we show that the latter possibility is rather exceptional.
More generally, we show that the set of the points where boundary conditions are lost,
can be prescribed to be arbitrarily close to any given open subset of the boundary.

\vskip 0.2cm

{\pc Keywords:}\ Diffusive Hamilton-Jacobi equation, viscosity solution, gradient blow-up, loss of boundary conditions

}}

\hskip 2mm \hbox{\box1}

\bigskip

{\bf 1. Introduction.}
\medskip

We consider the initial-boundary value problem for the diffusive Hamilton-Jacobi equation:
$$
\left\{\eqalign{
u_t-\Delta u&=|\nabla u|^p,\quad x\in\Omega,\ t>0,\cr
u(x,t)&=0,\quad x\in\partial\Omega,\ t>0,\cr
u(x,0)&=u_0(x),\quad x\in\Omega.
}\right.
\leqno(1.1)$$

Thoughout this article, we assume that $\Omega$ is a $C^{2+\alpha}$ smooth bounded domain of $\R^n$,
$p>2$ and 
$$\hbox{$u_0\in X:=\bigl\{v\in C^1(\overline\Omega)$; $v\ge 0\,$ and $v=0$ on $\partial\Omega\bigr\},$}$$
endowed with the $C^1$ norm.
This problem has been studied by many authors in the past decades (see e.g. [\QS, Chapter 40] and the references therein).
\smallskip

By standard theory [\Fri], it is known that problem (1.1) admits a maximal {\bf classical $C^1$ solution}
$u\ge 0$, such that $u\in C^{2,1}(\overline\Omega\times (0,T^*))$ and $u,\nabla u\in  C(\overline\Omega\times [0,T^*))$.
Here $T^*=T^*(u_0)\in (0,\infty]$ denotes the maximal existence time
and the differential equation and the boundary conditions are satisfied in the pointwise sense for $t\in(0,T^*)$.
Moreover, the solution satisfies the maximum principle estimate
$$\|u(t)\|_\infty\le \|u_0\|_\infty,\quad 0<t<T^*,$$
and the classical $C^1$ solution can only cease to exist through {\bf gradient blowup}:
$$T^*<\infty\quad\Longrightarrow\quad \lim_{t\to T^*}\|\nabla u(t)\|_\infty=\infty.$$
Actually, $\nabla u$ remains bounded away from the boundary and gradient blowup occurs only on $\partial\Omega$
(see [\SZ]). Furthermore it is known (see, e.g., [\Alaa, \ABG, \Sou]) that  $T^*<\infty$ whenever the initial data is suitably large.
We also recall that this phenomenon does not occur when $1\le p\le 2$.

\smallskip

On the other hand, it was proved in [\BaD] that problem (1.1) admits a unique {\bf global viscosity solution}
$u\in C(\overline\Omega\times [0,\infty))$,
where the boundary conditions have to be understood in the viscosity sense.
Throughout this article, we shall denote this solution by $u$ without risk of confusion, since the two solutions
coincide on $[0,T^*)$.
(The result in [\BaD] is actually valid for any $u_0\in C_0(\overline\Omega)$, but this need not concern us here.)
Moreover, $u$ is actually smooth away from the boundary, namely
$$u\in C^{2,1}(\Omega\times (0,\infty))$$
and it solves the PDE in (1.1) in the classical sense in $\Omega\times (0,\infty)$ (see Section~3 for details).
It was next proved in [\PZ] that for $t>T_0=T_0(\|u_0\|_\infty)$ sufficiently large,
$u$ is actually a classical solution again, namely $u(t)\in C^{2,1}(\overline\Omega\times (T_0,\infty))$ 
with $u(\cdot,t)=0$ on~$\partial\Omega$.

\smallskip

When gradient blowup occurs, the question of possible loss of boundary conditions for $t\ge T^*$ (hence actually in $[T^*,T_0]$)
has remained essentially open.
Namely, it is unknown whether or not $u$
satisfies the boundary conditions $u=0$ on $\partial\Omega\times[T^*,T_0]$ in the classical sense.
In what follows, we say that {\bf loss of boundary conditions} occurs at a point $x_0\in \partial\Omega$ if 
$u(x_0,t)>0$ for some $t\ge T^*$.

\smallskip

The goal of this article is to give some answers to this question.
A main conclusion is that loss of boundary conditions after gradient blowup
 {\bf may or may not occur, depending on the initial data.}
Furthermore, in case it occurs, the structure and size of the set of the points where boundary conditions are lost,
strongly depends on the initial data.
This is somewhat surprising and shows that the problem reveals a rather rich variety of phenomena.
\smallskip

Throughout this paper, we denote by $\varphi_1$ the first Dirichlet eigenfunction of $-\Delta$ in $\Omega$,
normalized by $\int_\Omega\varphi_1\, dx=1$.
\bigskip

\goodbreak
{\bf 2. Main results.}
\medskip

For any $u_0\in X$, we define the {\bf loss of boundary conditions set} by
$${\cal L}(u_0)=\bigl\{x_0\in \partial\Omega,\ u(x_0,t)>0 \hbox{ for some } t>0\bigr\}.$$
\smallskip

Our first result shows that there exist initial data for which
the loss of boundary conditions occurs {\bf everywhere on $\partial\Omega$},
and moreover can be achieved at a common time.

\medskip

\proclaim Theorem 1.
Let $p>2$. There exists $u_0\in X$ such that ${\cal L}(u_0)=\partial\Omega$
and that, moreover, 
$$u(x,t_0)\ge c_0,\quad x\in\partial\Omega,$$
for some $t_0,c_0>0$.
Furthermore, the same remains true for any $v_0\in X$ with $v_0\ge u_0$.

\medskip
Our next result shows that, at the opposite, there are gradient blowup solutions for which 
{\bf no loss of boundary conditions ever occurs.}
\medskip

\proclaim Theorem 2. 
Let $p>2$. There exists $u_0\in X$ such that $T^*(u_0)<\infty$ and ${\cal L}(u_0)=\emptyset$, i.e.,
$$u=0\quad\hbox{ on $\partial\Omega\times (0,\infty)$.}$$

At least in one space dimension, one can show that such $u_0$ are rather exceptional
(see Remark 2.2 for more comments).

\smallskip

\proclaim Theorem 3. 
Let $p>2$, $n=1$ and let $u_0$ be as in Theorem 2.
Let $v_0\in X$, with $v_0\not\equiv u_0$ and denote by $v$ the corresponding
global viscosity solution of (1.1).
\smallskip
(i) If $v_0\ge u_0$, then ${\cal L}(v_0)\neq\emptyset$. 
\smallskip
(ii) If $v_0\le u_0$, then $T^*(v_0)=\infty$.

\goodbreak

\medskip
Our last two results are concerned with some ''intermediate'' ranges of $u_0$. 
We first consider initial data which are large in an integral sense (hence need not be the same as those in Theorem 1) 
and show that, as the size grows larger, the loss of boundary conditions occurs "near" every point of $\partial\Omega$.
\medskip

\proclaim Theorem 4.
Let $p>2$. For any $\eps>0$, there exists a constant $M=M(\Omega,p,\eps)>0$ such that if $\int_\Omega u_0\varphi_1\, dx\ge M$,
then for any $x_0\in \partial\Omega$, we have ${\cal L}(u_0)\cap B_\eps(x_0)\neq\emptyset$.

\medskip

Notice that, due to the continuity of $u$ up to the boundary, ${\cal L}(u_0)$
is a (relatively) open subset of $\partial\Omega$.
Our last result shows that one can find solutions for which the loss of boundary conditions
occurs essentially only on a prescribed open subset of~$\partial\Omega$,
and at a common time.

\medskip
\proclaim Theorem 5.
Let $p>2$. Let $\omega$ be any open set of $\R^n$. Let $\eps>0$ and set $\omega_\eps=\omega+B_\eps(0)$.
There exists $u_0\in X$ such that
$$\omega\cap\partial\Omega \subset{\cal L}(u_0) \subset \omega_\eps\cap\partial\Omega$$
and that, moreover,
$$u(x,t_0)\ge c_0,\quad x\in\omega\cap\partial\Omega,
\leqno(\ConclLossCommon)$$
for some $t_0,c_0>0$.

\goodbreak

\medskip
{\bf Remark 2.1.} We note that some solutions with {\bf single-point} gradient blowup on the boundary (at $T^*(u_0)$)
may develop loss of boundary conditions on some {\bf open subset} of the boundary after $T^*(u_0)$.
Indeed, a closer inspection of the proof of Theorem~5 shows that one can construct $u_0$
which satisfy the conclusions of Theorem 5 
and at the same time verify the assumptions of [\LS, Theorem 1.1],
guaranteeing single-point gradient blowup on the boundary, for suitable domains of $\R^2$.

\medskip
{\bf Remark 2.2.} As shown by Theorem 3, the solutions constructed in Theorem 2 constitute (strong) {\bf thresholds}, 
realizing the transition from global classical existence to loss of boundary conditions. This parallels the phenomenon of transition from global existence to (complete) blowup for the reaction-diffusion equation 
$$u_t-\Delta u=u^p
\leqno(\RDeqn)$$
(see [\NST, \QS], and the recent work [\Qu] where the notion of {\it strong} threshold is studied). 
In this respect, gradient blowup without loss of boundary conditions plays the same role as
``incomplete blowup'' in the case of equation (\RDeqn), 
which is the threshold behavior for supercritical $p$ (i.e., $n\ge 3$ and $p>(n+2)/(n-2)$).

\medskip

For related results on the continuation of solutions after gradient blow-up, see [\FL, \FTW, \TT]. 
We refer to [\ABG, \FL, \CG, \ARS, \HM, \SV, \SZ, \GH, \LS, \PS]
for other aspects of gradient blowup phenomena, and to [\KPZ, \KS] for some physical background.

\goodbreak

\medskip
{\bf 3. Preliminaries}
\medskip

We set $Q=\Omega\times (0,\infty)$ and denote the function distance to the boundary by
$$\delta(x)={\rm dist}(x,\partial\Omega).$$
We set $\Omega_\eta=\{x\in\Omega;\ \delta(x)>\eta\}$ 
and recall that $\Omega_\eta$ is smooth for $\eta>0$ small.
Moreover, denoting by $\nu_\eta$ the outer normal unit vector and $d\sigma_\eta$ the surface measure on 
$\partial\Omega_\eta$, we have the property
$$\lim_{\eta\to 0}\int_{\partial\Omega_\eta} V\cdot\nu_\eta\,d\sigma_\eta=\int_{\partial\Omega} V\cdot\nu\,d\sigma,
\qquad V\in (C(\overline Q))^n.
\leqno(\ConvBoundary)$$
We shall also need the following uniform version of the Poincar\'e inequality (see, e.g., [\Ma]). 
Let $k\in [1,\infty)$. For each $\eps>0$, there exists
a constant $C=C(\Omega,\eps,k)>0$ such that
$$\int_\Omega |v|^k\le C \int_\Omega |\nabla v|^k\,dx,\qquad
v\in \displaystyle\bigcup_{x_0\in \partial\Omega} \{v\in W^{1,k}(\Omega);\ v_{|\partial\Omega\cap B_\eps(x_0)}=0\}
\leqno(\PoincareIneq)$$
($v=0$ being understood in the sense of traces if $v$ is not continuous).
\smallskip

We now turn to properties of the unique global viscosity solution $u$ of (1.1).
We refer to [\CIL, \BaD] for more details about viscosity solutions theory.
We first recall that the global viscosity solution can also be viewed as the limit of global classical solutions
of regularized problems.
Namely, for each integer $j\ge 1$, we set 
$$F_j(\xi)=\min\bigl(|\xi|^p,j^{p-2}|\xi|^2\bigr),\quad \xi\in \R^n,$$
and, for $u_0\in X$, 
 consider the problem 
$$
\left\{\eqalign{
\partial_tu_j-\Delta u_j&=F_j(\nabla u_j),\quad x\in\Omega,\ t>0,\cr
u_j(x,t)&=0,\quad x\in\partial\Omega,\ t>0,\cr
u_j(x,0)&=u_0(x),\quad x\in\Omega.
}\right.
\leqno(\approxpbm)
$$
Since each $F_j$ has at most quadratic growth, problem (\approxpbm) admits a unique global classical solution $u_j\ge 0$.
Moreover $u_j$ is nondecreasing with respect to $j$ by the comparison principle, 
and it is known (see [\CIL] and [\PZ]) 
that 
$$\lim_{j\to\infty}u_j(x,t)=u(x,t),\qquad (x,t)\in Q.$$

As a consequence of this approximation procedure, one for instance easily recovers the 
maximum principle estimate
$$\|u(t)-v(t)\|_\infty\le \|u_0-v_0\|_\infty,\quad t>0
\leqno(\contdep)$$
for all $u_0, v_0\in X$ 
(which yields in particular the continuous dependence in $L^\infty$).

Next, as a consequence of uniform interior estimates for the approximating solutions $u_j$, one shows that
$$u\in C^{2,1}(Q)
\leqno(\classreg)$$
and that $u$ solves the PDE in (1.1) in the classical sense in $Q$.
For that purpose, by standard parabolic regularity, it suffices to prove that $\nabla u_j$ is bounded on compact subsets 
of $Q$, independently of $j$. 
Such a bound can be proved by a Bernstein argument with cut-off
(see e.g. [\Lio] in the elliptic case and [\SZ] in the parabolic case; more specifically, 
this follows from a simple modification of the proof of [\SZ, Theorem 3.2]).

\medskip

Moreover, we have the following time-derivative estimate.

\proclaim Lemma 3.1. Let $u_0\in X$ and let $u$ be the corresponding global viscosity solution of~(1.1).
Then, for all $t>0$ we have $u_t(\cdot,t)\in L^\infty(\Omega)$.
Moreover,  for all $t_0>0$, there exists a constant $C(t_0)>0$ such that
$$\|u_t(t)\|_\infty\le C(t_0),\quad t\ge t_0.
\leqno(\estimut)$$

{\it Proof.} We may assume without loss of generality that $t_0\in (0,T^*(u_0))$.
Let $t_0\le t<t+h<T^*(u_0))$.
By estimate (\contdep), we have
$$\|u(t+h)-u(t)\|_\infty\le \|u(t_0+h)-u(t_0)\|_\infty.$$
Recall that $u\in C^{2,1}(\overline\Omega\times (0,T^*))\cap C^{2,1}(\Omega\times (0,\infty))$.
Dividing by $h$ and letting $h\to 0$, we deduce that
$$\|u_t(t)\|_\infty\le \|u_t(t_0)\|_\infty,$$
and the lemma is proved. \fin

On the other hand, we know that $\nabla u$ also satisfies the following Bernstein estimate: for each $\tau>0$, there exists a constant $C(\tau)>0$ such that
$$|\nabla u(x,t)|\le C(\tau)\delta^{-1/(p-1)}(x),\quad x\in \Omega,\ t\ge \tau.
\leqno(\Bernstein)$$
This is proved in [\SZ] for classical solutions, i.e., on $(0,T^*(u_0))$, but the proof remains valid for the global viscosity solution,
using (\classreg) along with estimate (\estimut).

\smallskip

Finally, we give the following lemma, which will be useful for the proof of Theorem~3.

\proclaim Lemma 3.2. Let $u_0, v_0\in X$ and $\lambda>1$ be such that $v_0\ge \lambda u_0$ 
and denote by $u,v$ the corresponding global viscosity solutions of (1.1).
Then 
$$v\ge \lambda u\quad\hbox{ in $\overline\Omega\times (0,\infty)$.}$$

{\it Proof.} Let $j\ge 1$ and let $u_j, v_j$ be the solutions of the approximating problems (\approxpbm).
Setting $\underline u_j=\lambda u_j$, we see that
$$\eqalign{
\partial_t\underline u_j-\Delta \underline u_j-F_j(\nabla \underline u_j)
&=\lambda\Bigl[\min\bigl(|\nabla u_j|^p,j^{p-2}|\nabla u_j|^2\bigr)
-\min\bigl(\lambda^{p-1}|\nabla u_j|^p,j^{p-2}\lambda|\nabla u_j|^2\bigr)\Bigr]\cr
&\le\lambda\Bigl[\min\bigl(|\nabla u_j|^p,j^{p-2}|\nabla u_j|^2\bigr)
-\min\bigl(|\nabla u_j|^p,j^{p-2}|\nabla u_j|^2\bigr)\Bigr] \cr
&=0=\partial_tv_j-\Delta v_j-F_j(\nabla v_j)
}$$
in $Q$. We deduce from the comparison principle that $\underline u_j\le v_j$ in $Q$
and the result follows by passing to the limit $j\to\infty$. \fin

\medskip
{\bf 4. Proof of Theorem 4.}
\medskip

We first prove Theorem 4, since the result is (independently) used in the proof of Theorem~1.
We adapt eigenfunction arguments used in [\Alaa, \Sou] to prove gradient blowup for weak or classical solutions. 
It turns out that these arguments can be modified to establish the loss of boundary conditions for global viscosity solutions,
making use of the Bernstein estimate (\Bernstein).

Recall that we denote by $\varphi_1$ the first Dirichlet eigenfunction of $-\Delta$ in $\Omega$,
normalized by $\int_\Omega\varphi_1\, dx=1$ and let $\lambda_1>0$ be the corresponding eigenvalue.
Let $0<\tau<t<\infty$ and let $\eta>0$ small. 
Since we only have $u\in C^{2,1}(Q)\cap C(\overline Q)$, we cannot directly integrate in $\Omega$.
Instead, we multiply the PDE in (1.1) by $\varphi_1$ and integrate by parts
on $\Omega_\eta$. This yields
$$
\eqalign{
\Bigl[\int_{\Omega_\eta} u\varphi_1\,dx\Bigr]_\tau^t
&=\int_\tau^t\int_{\Omega_\eta} \varphi_1\Delta u\,dxds+\int_\tau^t\int_{\Omega_\eta} |\nabla u|^p\varphi_1\,dxds \cr
&=\int_\tau^t\int_{\Omega_\eta} u\Delta  \varphi_1\,dxds+\int_\tau^t\int_{\partial\Omega_\eta} (\varphi_1\nabla u-u\nabla\varphi_1)\cdot\nu_\eta\,d\sigma_\eta ds\cr
&\qquad\qquad\qquad +\int_\tau^t\int_{\Omega_\eta} |\nabla u|^p\varphi_1\,dxds. \cr
}$$
Recall that
$$c_1\delta(x)\le\varphi_1(x)\le c_2\delta(x),\qquad x\in\Omega,
\leqno(\propdeltaA)$$
and
$$\int_\Omega \delta^{-\beta}(x)\,dx<\infty,\quad\hbox{ for all $\beta\in (0,1)$}
\leqno(\propdeltaB)$$
(see e.g. [\Sou]).
Using (\Bernstein), (\propdeltaA) and (\ConvBoundary), we obtain
$$
\eqalign{\Bigl|\int_\tau^t\int_{\partial\Omega_\eta} \varphi_1\nabla u\cdot\nu_\eta\,d\sigma_\eta ds\Bigr |
&\le C(\tau)t\int_{\partial\Omega_\eta}\delta^{(p-2)/(p-1)}(x)\,d\sigma_\eta \cr
&\le C(\tau)t\eta^{(p-2)/(p-1)}\int_{\partial\Omega_\eta}\,d\sigma_\eta \cr
& \le C(\tau)t\eta^{(p-2)/(p-1)}\to 0,\quad\hbox{as $\eta\to 0.$}
}$$
Also we note that, for all $t>0$, we have
$$\int_\Omega |\nabla u(t)|^p\varphi_1\,dx\le C(t)\int_\Omega \delta^{-p/(p-1)}(x)\delta(x)\,dx
=C(t)\int_\Omega \delta^{-1/(p-1)}(x)\,dx<\infty,
\leqno(\integralvarphiA)$$
owing to (\Bernstein), (\propdeltaA) and (\propdeltaB).
Using (\ConvBoundary) and the fact that $u\in C(\overline Q)$, we may pass to the limit $\eta\to 0$ to get
$$
\eqalign{
\Bigl[\int_\Omega u\varphi_1\,dx\Bigr]_\tau^t
&=-\lambda_1\int_\tau^t\int_\Omega u \varphi_1\,dxds \cr
&\qquad -\int_\tau^t\int_{\partial\Omega} u\partial_\nu\varphi_1\,d\sigma ds+\int_\tau^t\int_\Omega |\nabla u|^p\varphi_1\,dxds. \cr
}$$
Using $u\ge 0$ and $\partial_\nu\varphi_1\le 0$ on $\partial\Omega$,
and then passing to the limit $\tau\to 0$, we get, for all $t>0$,
$$\int_\Omega u(t)\varphi_1\,dx\ge \int_\Omega u_0\varphi_1\,dx
+\int_0^t\int_\Omega |\nabla u|^p\varphi_1\,dxds
-\lambda_1\int_0^t\int_\Omega u \varphi_1\,dxds,
\leqno(\integralvarphi)
$$
hence in particular the finiteness of the integral of the gradient term in (\integralvarphi).
Let $k\in [1,p/2)$. By H\"older's inequality, we have
$$
\int_\Omega |\nabla u|^k\,dx=\int_\Omega |\nabla u|^k\varphi_1^{k/p}\varphi_1^{-k/p}\,dx
\le \Bigl(\int_\Omega |\nabla u|^p\varphi_1\,dx\Bigr)^{k/p} \Bigl(\int_\Omega \varphi_1^{-k/(p-k)}\,dx\Bigr)^{(p-k)/p},
$$
hence
$$
\Bigl(\int_\Omega |\nabla u|^k\,dx\Bigr)^{p/k}\le C(k)\int_\Omega |\nabla u|^p\varphi_1\,dx.
\leqno(\integralvarphiB)$$
owing to (\propdeltaA) and (\propdeltaB). In particular, in view of (\integralvarphiA) and of $u\in C(\overline Q)$, we have
$$u(t)\in W^{1,k}(\Omega),\quad\hbox{ for all $t>0$.}
\leqno(\uWonep)$$

Now assume that there exist $\eps>0$ and $x_0\in \partial\Omega$ such that
$u=0$ on $(\partial\Omega\cap B_\eps(x_0))\times (0,\infty)$. 
Fixing any $k\in (1,p/2)$, and taking (\uWonep) into account, we may
apply the Poincar\'e inequality (\PoincareIneq). This along with H\"older's inequality and (\integralvarphiB) yields
$$
\Bigl(\int_\Omega u\varphi_1\,dx\Bigr)^p
\le \Bigl(\int_\Omega u^k\,dx\Bigr)^{p/k}
\le C(\eps)\Bigl(\int_\Omega |\nabla u|^k\,dx\Bigr)^{p/k}
\le C(\eps)\int_\Omega |\nabla u|^p\varphi_1\,dx.
$$
Going back to (\integralvarphi), it follows that, for all $t>0$,
$$\int_\Omega u(t)\varphi_1\,dx
\ge \int_\Omega u_0\varphi_1\,dx
+c_0\int_0^t \biggl[\Bigl(\int_\Omega u\varphi_1\,dx\Bigr)^p-c_1^p\biggr]\,ds,
\leqno(\integralvarphiD)$$
for some constants $c_0, c_1>0$ depending on $\eps$.
Assume that
$$\int_\Omega u_0\varphi_1\,dx\ge 2c_1.$$
It then easily follows from (\integralvarphiD) that 
$\int_\Omega u(t)\varphi_1\,dx\ge \int_\Omega u_0\varphi_1\,dx\ge 2c_1$ for all $t>0$.
Consequently
$$\int_\Omega u(t)\varphi_1\,dx
\ge \int_\Omega u_0\varphi_1\,dx
+{c_0\over 2}\int_0^t \Bigl(\int_\Omega u\varphi_1\,dx\Bigr)^p\,ds=:H(t),\quad t>0,
$$
hence $H'(t)\ge (c_0/2)H^p$, which implies the finite time blowup of $\int_\Omega u(t)\varphi_1\,dx$.
But this is a contradiction with the fact that $u\in C(\overline Q)$ 
(or with the estimate $\int_\Omega u(t)\varphi_1\,dx\le \|u(t)\|_\infty\le \|u_0\|_\infty$).
\fin
\bigskip
\goodbreak

\medskip
{\bf 5. Proof of Theorem 1.}
\medskip

We shall modify an argument from [\LS] based on a radial auxiliary problem and a scaling argument.
Consider the auxiliary problem
$$
\left\{\eqalign{
V_t-\Delta V&=|\nabla V|^p,\quad x\in B_1(0),\ t>0,\cr
V(x,t)&=0,\quad x\in \partial B_1(0),\ t>0,\cr
V(x,0)&=V_0(x),\quad x\in B_1(0),
}\right.
\leqno(\AuxilPbmV)
$$
where $V_0\in C^1(\overline B_1(0))$, with $V_0$ radially symmetric and supported in $B_{1/2}(0)$.
As a consequence of Theorem~4, proved in the previous section, we may choose $V_0$ such that loss of boundary conditions occurs for $V$.
Since $V$ is radially symmetric, it follows that there exist $t_0, c_0>0$ such that
$$V(x,t_0)=c_0\quad\hbox{ for all $x\in \partial B_1(0)$.}$$

Next, since $\partial\Omega$ is smooth, one can find $\rho>0$ such that,
for all $x_0\in \partial\Omega$, there exists $x_1=x_1(x_0)$ such that 
$$B_\rho(x_1)\subset\Omega\quad\hbox{ and }\quad \partial\Omega\cap\partial B_\rho(x_1)=\{x_0\}.
\leqno(\TangentBall)$$
We now use the scale invariance of the equation and set 
$$w(x_0;x,t)=\rho^\beta V\bigl(\rho^{-1}(x-x_1),\rho^{-2}t\bigr),\quad (x,t)\in \overline B_\rho(x_1)\times [0,\infty),$$
with $\beta=(p-2)/(p-1)$. 
A straightforward computation shows that $w(x_0;\cdot,\cdot)$ is the solution of (\AuxilPbmV)
with $B_1(0)$ replaced by $B_\rho(x_1)$ and $V_0$ replaced by $\rho^\beta V_0(\rho^{-1}(x-x_1))$.
Now choose $u_0\in X$ such that 
$$u_0\ge \rho^\beta\|V_0\|_\infty\quad\hbox{ in $\Omega'_\rho:=\{x\in\Omega;\ \rho/2\le \delta(x)\le 3\rho/2\}$.}$$
For any $x_0\in \partial\Omega$, the function $w(x_0;\cdot,0)$ is supported in $B_{\rho/2}(x_1)\subset\Omega'_\rho$,
owing to (\TangentBall), hence $u_0\ge w(x_0;\cdot,0)$ in $\overline B_\rho(x_1)$.
By the comparison principle, it follows that $u\ge w(x_0;\cdot,\cdot)$ in $\overline B_\rho(x_1)\times [0,\infty)$,
hence in particular, 
$$u(x_0,\rho^2t_0)\ge w(x_0;x_0,\rho^2 t_0)=\rho^\beta V\bigl(\rho^{-1}(x_0-x_1),t_0\bigr)=c_0\rho^\beta>0.$$
The conclusion for $u_0$ follows.
The assertion for $v_0\ge u_0$ is then an immediate consequence of the comparison principle. \fin

\medskip
{\bf 6. Proof of Theorem 2.}
Fix $\phi\in X$, $\phi\not\equiv 0$ and, for $\lambda>0$, consider $u_\lambda$ the solution of (1.1)
with initial data $u_0=\lambda\phi$.
By, e.g., [\Sou] we know that $T^*(\lambda\phi)=\infty$ for $\lambda$ small and $T^*(\lambda\phi)<\infty$ for $\lambda$ large.
We may thus define
$$\lambda^*=\inf\{\lambda>0;\ T^*(\lambda\phi)<\infty\}\in (0,\infty).$$
We shall prove that $u_{\lambda^*}$ has the desired properties.

First, since $u_\lambda=0$ on $\partial\Omega\times (0,\infty)$ for all $\lambda\in(0,\lambda^*)$,
it follows from the $L^\infty$ continuous dependence estimate (\contdep) that
$$u_{\lambda^*}=0\quad\hbox{ on $\partial\Omega\times (0,\infty)$.}$$

Next, by [\Sou], the trivial solution is asymptotically stable in $X$.
Namely, there exists $\eps_0=\eps_0(\Omega,p)>0$ such that,
for any $v_0\in X$, 
$$\|v_0\|_X\le \eps_0\quad\Longrightarrow\quad T^*(v_0)=\infty\ \hbox{ and }\ \lim_{t\to\infty}\|v(t)\|_X=0.$$
On the other hand, by [\PZ], we know that there exists $t_0>0$ such that 
$$u_{\lambda^*}(t)\in X\ \hbox{ for all $t\ge t_0\ \ $ and }\ \lim_{t\to\infty} \|u_{\lambda^*}(t)\|_X=0.$$
Consequently, there exists $t_1>t_0$ such that $\|u_{\lambda^*}(t_1)\|_X<\eps_0(\Omega,p)$.

Now assume for contradiction that $T^*(\lambda^*\phi)=\infty$.
Then by continuous dependence in $X$ of classical solutions, 
there exists $\eta>0$ such that
$$\hbox{ if $|\lambda-\lambda^*|<\eta$, then $T^*(\lambda\phi)>t_1$ and 
$\|u_{\lambda}(t_1)\|_X<\eps_0$.}$$
By the above asymptotically stability property,
it follows in particular that $T^*(\lambda\phi)=\infty$ for all $\lambda\in (\lambda^*,\lambda^*+\eta)$.
But this contradicts the definition of $\lambda^*$. The proof is complete.
\fin
\medskip

\medskip
{\bf 7. Proof of Theorem 3.}
Assume without loss of generality that $\Omega=(-1,1)$. Set
$$\beta=(p-2)/(p-1),\qquad c_p=(p-2)^{-1}(p-1)^{-1/(p-1)}.$$
For any $w_0\in X$, denoting by $w$ the corresponding solution of (1.1),
we know from [\CG] and [\QS, Theorem~40.14] that, if $T^*(w_0)<\infty$, then
$$\lim_{x\to x_0}{w(x,T^*(w_0))\over  \delta^{\beta}(x)}=c_p,\quad\hbox{ for some $x_0\in \{-1,1\}$}.
\leqno(\ratiowA)$$

Next, for any $t>0$ and any $x_0\in \{-1,1\}$, we claim that
$$w(x_0,t)=0\quad\Longrightarrow\quad \limsup_{x\to x_0}{w(x,t)\over  \delta^{\beta}(x)}\le c_p.
\leqno(\ratiowB)$$
Consider the case $x_0=-1$ (the other case being similar). For a fixed $t>0$, 
we let 
$$y(x)=(w_x(x,t)-C_1 (x+1))_+,$$
where $C_1=C(t)$ is given by (\estimut).
The function $y$ satisfies
$$y'+y^p=(w_{xx}-C_1)\chi_{\{w_x>C_1 (x+1)\}}+(w_x-C_1(x+1))_+^p,
\quad\hbox{for a.e. $x\in (-1,0]$.}$$
For each $x$ such that $w_x(x,t)>C_1(x+1)$, 
we have 
$$(y'+y^p)(x)\leq (w_{xx}-C_1+|w_x|^p)(x,t)\le 0$$ 
by (1.1) and (\estimut). 
Therefore, we have $y'+y^p\le 0$ a.e. on $(-1,0]$. By integration, it follows that
$y(x)\leq ((p-1)(x+1))^{-{1\over p-1}}$, hence
$w_x(x,t)\leq ((p-1)(x+1))^{-{1\over p-1}}+C_1$ on $(-1,0]$.
Assuming $w(-1,t)=0$, a further integration then yields
$$w(x,t)\leq c_p(x+1)^\beta+C_1(x+1),\qquad x\in(-1,0],$$
and claim (\ratiowB) is proved.

Let now $u_0\in X$ be such that $T^*(u_0)<\infty$ and 
$$u=0\quad\hbox{ on $\partial\Omega\times (0,\infty).$}
\leqno(\HypNoLoss)$$
We first prove assertion (i) and consider $v_0\in X$ such that $v_0\ge u_0$ and $v_0\not\equiv u_0$.
Pick $t_0$ such that $0<t_0<T^*(v_0)\le T^*(u_0)$.
It follows easily from the Hopf Lemma that $v(\cdot,t_0)\ge \lambda u(\cdot,t_0)$ in $\Omega$
for some $\lambda>1$. By Lemma~3.2, we deduce that $v\ge \lambda u$ in $\overline\Omega\times [t_0,\infty)$.
Next applying~(\ratiowA) with $w=u$, it follows that there exists $x_0\in \{-1,1\}$ such that
$$\limsup_{x\to x_0}{v(x,T^*(u_0))\over  \delta^{\beta}(x)}\ge \lambda c_p>c_p.$$
As a consequence of (\ratiowB) applied with $w=v$, we deduce that $v(x_0,T^*(u_0))>0$.

To prove assertion (ii), we consider $v_0\in X$ such that $v_0\le u_0$ and $v_0\not\equiv u_0$.
Arguing similarly as before, we deduce that $v\le \lambda u$ in $\overline\Omega\times [t_0,\infty)$ for some $\lambda<1$.
By (\HypNoLoss) and (\ratiowB) applied with $w=u$, for any $t>0$ and any $x_0\in \{-1,1\}$, it follows that
$$\limsup_{x\to x_0}{v(x,t)\over  \delta^{\beta}(x)}\le \lambda c_p<c_p.$$
As a consequence of (\ratiowA) applied with $w=v$, we deduce that $T^*(v_0)=\infty$.
The result is proved. \fin

\medskip
{\bf 8. Proof of Theorem 5.}
First, following [\LS, Lemma 2.3], we fix a smooth function $h\geq 0$ in $\R^n$ such that
$$
    h(x)=\cases{
                  1, & $x\in \omega$ \cr
                  \noalign{\vskip 1mm}
                  0, & $x\in\R^n\setminus \omega_\eps$}
\leqno(\PsiElliptA)$$
and consider the elliptic problem
$$
\left\{\eqalign{
  -\Delta \psi &=1,\quad x\in\Omega,\cr
  \psi& =h,\quad x\in\partial\Omega.\cr
}\right.
\leqno(\PsiElliptB)$$
We have
$$
    -\Delta(c_1\psi)=c_1\geq |\nabla(c_1\psi)|^p,\
    \ \ x\in\Omega,
\leqno(\PsiEllipt)$$
with $c_1:=\|\nabla\psi\|_{L^\infty(\Omega)}^{-p/(p-1)}>0$. 

\goodbreak

Next, by the continuity of $\psi$ in $\overline{\Omega}$, we may find $\rho\in (0,\eps/3)$ such that
$$\psi\ge 1/2\quad\hbox{ in $\{x\in \overline\Omega;\ {\rm dist}(x,\omega\cap\partial\Omega)\le 2\rho\}$}.
\leqno(\PsiLarge)$$
Taking $\rho$ smaller and owing to the regularity of $\partial\Omega$, we may also assume
that for all $x_0\in \partial\Omega$, there exists a point $x_1=x_1(x_0)$ such that 
$$B_\rho(x_1)\subset\Omega\quad\hbox{ and }\quad \partial\Omega\cap\partial B_\rho(x_1)=\{x_0\}.
\leqno(\TangentBallB)$$
Now let $V_0$ be as in the proof of Theorem~1. Taking $\rho$ even smaller, we may also assume that 
$$\rho^\beta\|V_0\|_\infty<{c_1\over 2}.$$
We may thus choose $u_0\in X$ such that 
$$u_0=0\quad\hbox{ in $\{x\in \overline\Omega;\ {\rm dist}(x,\omega\cap\partial\Omega)\ge2\rho\},$}$$
$$u_0\le {c_1\over 2}\quad\hbox{ in $\{x\in \overline\Omega;\ {\rm dist}(x,\omega\cap\partial\Omega)<2\rho\}$}$$
and
$$u_0\ge \rho^\beta\|V_0\|_\infty
\quad\hbox{ in $\Omega''_\rho:=\{x\in\Omega;\ \delta(x)\ge \rho/2\ $ and $\ {\rm dist}(x,\omega\cap\partial\Omega)\le 3\rho/2\}$.}
\leqno(\condVzero)$$
In particular, in view of (\PsiLarge), we have
$$u_0\le c_1\psi \quad\hbox{ in $\overline\Omega$}.$$
By (\PsiEllipt) and the comparison principle, we deduce that
$$u\le c_1\psi \quad\hbox{ in $\overline\Omega\times (0,\infty)$},$$
hence in particular ${\cal L}(u_0) \subset \omega_\eps\cap\partial\Omega$, due to (\PsiElliptA), (\PsiElliptB).

On the other hand, for each $x_0\in \omega\cap\partial\Omega$,
we can prove the loss of boundary conditions at $x_0$ by arguing similarly as in the proof of Theorem~1.
Namely, recalling (\TangentBallB), we set
$$w(x_0;x,t)=\rho^\beta V\bigl(\rho^{-1}(x-x_1),\rho^{-2}t\bigr),\qquad (x,t)\in \overline B_\rho(x_1)\times [0,\infty),$$
where $V$ is as in the proof of Theorem~1 and $\beta=(p-2)/(p-1)$. 
The function $w(x_0;\cdot,0)$ is supported in $B_{\rho/2}(x_1)\subset\Omega''_\rho$,
owing to (\TangentBallB), hence $u_0\ge w(x_0;\cdot,0)$ in $\overline B_\rho(x_1)$ by (\condVzero).
By the comparison principle, it follows that $u\ge w(x_0;\cdot,\cdot)$ in $\overline B_\rho(x_1)\times [0,\infty)$,
hence in particular, 
$$u(x_0,\rho^2t_0)\ge w(x_0;x_0,\rho^2 t_0)=\rho^\beta V\bigl(\rho^{-1}(x_0-x_1),t_0\bigr)=c_0\rho^\beta>0.$$
Therefore, $\omega\cap\partial\Omega\subset {\cal L}(u_0)$
and (\ConclLossCommon) holds.
The theorem is proved. \fin

\bigskip 
\noindent{\bf Acknowledgement.} \ 
Most of this work was done during a visit of Ph.~Souplet at the 
Dipartimento di Matematica of Universit\`a di Roma Tor Vergata in April 2016.
He wishes to thank this institution for the kind hospitality.

\goodbreak

{\baselineskip=11pt \parindent=0.75cm
\font\pc=cmcsc9
\font\rmn=cmr9
\font\sln=cmsl9
\font\rmb=cmbx8 scaled 1125 \rm
\font\it=cmti9

\bigskip\bigskip
\centerline{\rmb REFERENCES}
\vskip 4mm

\rmn \eightpoint

\item{[\Alaa]} Alaa, N.,
             {\sln Weak solutions of quasilinear parabolic equations with measures as initial data},
             Ann. Math. Blaise Pascal~3 (1996), 1--15.
\smallskip
              
\item{[\ABG]} Alikakos, N.D., Bates, P.W., Grant, C.P.,
              {\sln Blow up for a diffusion-advection equation},
              Proc. Roy. Soc. Edinburgh Sect. A 113 (1989), 181--190.
\smallskip
              
\item{[\ARS]} Arrieta, J.M., Rodriguez-Bern\'al, A., Souplet, Ph.,
             {\sln Boundedness of global solutions for nonlinear parabolic equations
             involving gradient blow-up phenomena},
              Ann. Sc. Norm. Super. Pisa Cl. Sci. (5) 3 (2004), 1--15.
\smallskip
            
\item{[\BaD]} Barles, G., Da Lio, F.,
             {\sln On the generalized Dirichlet problem for viscous Hamilton-Jacobi equations},
             J. Math. Pures Appl. 83 (2004), 53--75.
\smallskip

\item{[\BD]} Benachour, S., Dabuleanu, S.,
              {\sln The mixed Cauchy-Dirichlet problem for a viscous Hamilton-Jacobi equation},
              Adv. Differential Equations 8 (2003), 1409--1452.
\smallskip
      
\item{[\CG]} Conner, G.R., Grant, C.P.,
              {\sln Asymptotics of blowup for a convection-diffusion equation with conservation},
            Differential Integral Equations 9 (1996), 719--728.
            
\smallskip

\item{[\CIL]} Crandall, M., Ishii, H., Lions, P.L.,
              {\sln User's guide to viscosity solutions of second order partial differential equation},
            Bull. Amer. Math. Soc. 27 (1992), 1--67.
            
\smallskip
             
\item{[\FL]} Fila, M., Lieberman, G.M.,
            {\sln Derivative blow-up and beyond for quasilinear parabolic equations},
             Differential Integral Equations 7 (1994), no. 3-4, 811--821.
\smallskip

\item{[\FTW]}
     		Fila, M., Taskinen, J., Winkler, M.,
            {\sln Convergence to a singular steady-state of a parabolic equation with gradient blow-up},
		Appl. Math. Letters 20 (2007) 578--582.

\smallskip
    
\item{[\Fri]} Friedman, A., Partial differential equations of parabolic type,
             Prentice-Hall, 1964.
       
\smallskip

\item{[\GH]} Guo, J.-S., Hu, B.,
              {\sln Blowup rate estimates for the heat equation with a nonlinear gradient source term},
              Discrete Contin. Dyn. Syst. 20 (2008), 927--937.        
\smallskip
 
\item{[\HM]} Hesaaraki, M., Moameni, A.,
              {\sln Blow-up positive solutions for a family of nonlinear parabolic equations
             in general domain in $\R^N$},
             Michigan Math. J. 52 (2004), 375--389.           
\smallskip
 
\item{[\KPZ]} Kardar, M., Parisi, G., Zhang, Y.C.,
             {\sln Dynamic scaling of growing interfaces},
              Phys. Rev. Lett. 56 (1986), 889--892.
\smallskip
 
\item{[\KS]} Krug, J., Spohn, H.,
           {\sln Universality classes for deterministic surface growth},
            Physical Review A 38 (1988), 4271--4283.
\smallskip
   
\item{[\LS]} Li, Y.-X., Souplet, Ph.,
              {\sln Single-Point Gradient Blow-up on the Boundary for Diffusive Hamilton-Jacobi Equations in Planar Domains},
             Comm. Math. Phys. 293 (2009), 499--517.
\smallskip

\item{[\Lio]} Lions,~P.-L.
{\sln Quelques remarques sur les probl\`emes elliptiques quasilin\'eaires du second ordre.} (French)
J. Analyse Math. 45  (1985), 234--254.

\smallskip 
  
\item{[\Ma]} Maz'ya, V.G., Sobolev spaces. Springer 1985.

\smallskip 

\item{[\NST]} Ni, W.-M., Sacks, P., Tavantzis, J.,
{\sln On the asymptotic behavior of solutions of certain quasilinear parabolic equations},
J. Differential Equations 54 (1984), 97--120.

\smallskip 

\item{[\PS]} Porretta, A., Souplet, Ph.,
   {\sln The profile of boundary gradient blow-up for the diffusive Hamilton-Jacobi equation},
International Math. Research Notices (2016), in press, doi: 10.1093/imrn/rnw154.

\smallskip 

\item{[\PZ]} Porretta, A., Zuazua, E.,
   {\sln Null controllability of viscous Hamilton-Jacobi equations},
Ann. Inst. H. Poincar\'e Anal. Non Lin\'eaire 29 (2012), 301--333.

\smallskip 
  
\item{[\Qu]} Quittner, P., Threshold and strong threshold solutions of a semilinear parabolic equation, Preprint,
ArxiV 1605.07388 (2016).

\smallskip 

\item{[\QS]} Quittner, P., Souplet, Ph.,
          Superlinear parabolic problems. Blow-up, global existence and steady states,
              Birkh\"{a}user Advanced Texts: Basel Textbooks, Birkh\"{a}user Verlag, Basel, 2007.

\smallskip

\item{[\Sou]} Souplet, Ph.,
              {\sln Gradient blow-up for multidimensional nonlinear parabolic equations
              with general boundary conditions},
              Differential Integral Equations 15 (2002),  237--256.
\smallskip

\item{[\SV]} Souplet, Ph., V\'azquez, J.L.,
              {\sln Stabilization towards a singular steady state with gradient blow-up for a convection-diffusion problem},
Discrete Contin. Dyn. Syst. 14 (2006), 221--234.

\item{[\SZ]} Souplet, Ph., Zhang, Q.S.,
              {\sln Global solutions of inhomogeneous Hamilton-Jacobi equations},
              J. Anal. Math. 99 (2006), 355--396.
              
\smallskip

\item{[\TT]} Tabet Tchamba, Th.,
              {\sln Large time behavior of solutions of viscous Hamilton-Jacobi equations with superquadratic Hamiltonian},
               Asymptot. Anal. 66 (2010), 161--186. 
\bye